# SIZE, POWER AND FALSE DISCOVERY RATES

By Bradley Efron

*Stanford University*

Modern scientific technology has provided a new class of large-scale simultaneous inference problems, with thousands of hypothesis tests to consider at the same time. Microarrays epitomize this type of technology, but similar situations arise in proteomics, spectroscopy, imaging, and social science surveys. This paper uses false discovery rate methods to carry out both size and power calculations on large-scale problems. A simple empirical Bayes approach allows the false discovery rate (fdr) analysis to proceed with a minimum of frequentist or Bayesian modeling assumptions. Closed-form accuracy formulas are derived for estimated false discovery rates, and used to compare different methodologies: local or tail-area fdr's, theoretical, permutation, or empirical null hypothesis estimates. Two microarray data sets as well as simulations are used to evaluate the methodology, the power diagnostics showing why nonnull cases might easily fail to appear on a list of "significant" discoveries.

**1. Introduction.** Large-scale simultaneous hypothesis testing problems, with hundreds or thousands of cases considered together, have become a familiar feature of current-day statistical practice. Microarray methodology spearheaded the production of large-scale data sets, but other "high throughput" technologies are emerging, including time of flight spectroscopy, proteomic devices, flow cytometry and functional magnetic resonance imaging.

Benjamini and Hochberg's seminal paper [3] introduced false discovery rates (Fdr), a particularly useful new approach to simultaneous testing. Fdr theory relies on $p$-values, that is on null hypothesis tail areas, and as such operates as an extension of traditional frequentist hypothesis testing to simultaneous inference, whether involving just a few cases or several thousand. Large-scale situations, however, permit another approach: empirical Bayes methods can bring Bayesian ideas to bear without the need for strong

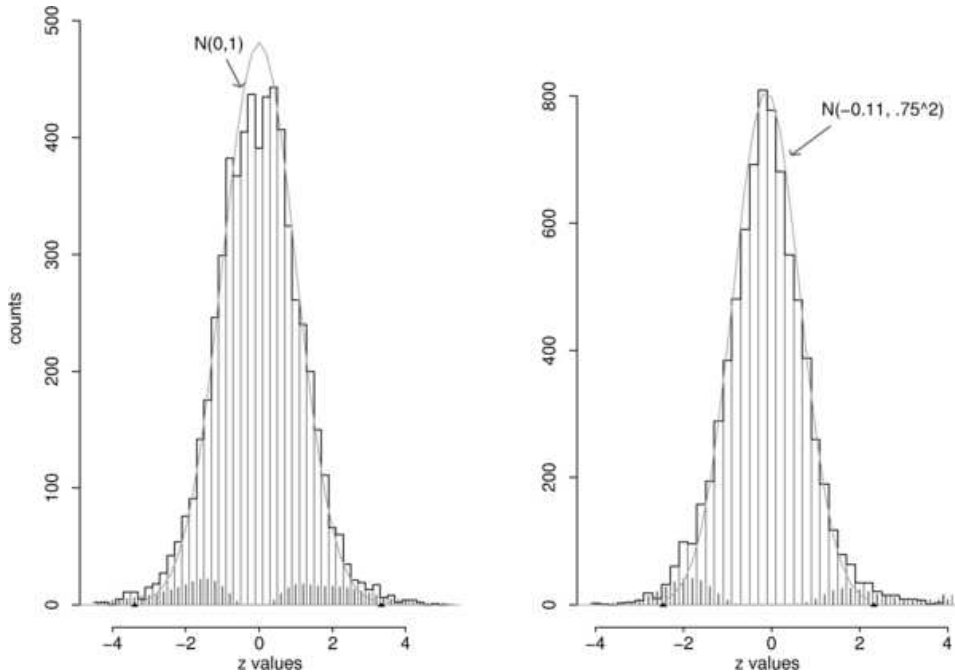

FIG. 1. *Histograms of z values from two microarray experiments. Left panel, prostate data, comparison of* 50 *nontumor subjects with* 52 *tumor patients for each of* 6033 *genes; Singh et al.* [31]. *Right panel, HIV data, comparison of* 4 *HIV negative subjects with* 4 *HIV positive patients for* 7680 *genes; van't Wout et al.* [34], *discussed in* [16]. *The central peak of the prostate data histogram closely follows the theoretical* $N(0,1)$ *null density (solid curve), but the HIV histogram is substantially too narrow. Short vertical bars are estimated nonnull counts, useful for power calculations, as discussed in Section* 3. *Estimated null proportion* $p_0$ *equals* 0.93 *in both experiments.*

Bayesian or frequentist assumptions. This paper pursues large-scale false discovery rate analysis from an empirical Bayes point of view, with the goal of providing a versatile methodology for both size and power considerations.

The left panel of Figure 1 concerns a microarray example we will use to introduce our main ideas. 102 microarrays, 50 from nontumor subjects and 52 from prostate cancer patients, each measured expression levels for the same $N = 6033$ genes. Each gene yielded a two-sample $t$-statistic $t_i$ comparing tumor versus nontumor men, which was then transformed to a $z$ value,

$$(1.1) \qquad z_i = \Phi^{-1}(F_{100}(t_i)),$$

where $F_{100}$ is the cumulative distribution function (c.d.f.) of a Student's $t$ distribution with 100 degrees of freedom, and $\Phi$ is the standard normal c.d.f.

We expect $z_i$ to have nearly a $N(0,1)$ distribution for "null" genes, the ones behaving similarly in tumor and nontumor situations. The left his-togram looks promising in this regard: its large central peak, which is nicely proportional to a $N(0,1)$ density, charts the presumably large majority of null genes, while the heavy tails suggest some interesting "nonnull" genes, those responding differently in the two situations, the kind the study was intended to detect.

NOTE. It is not necessary that the $z_i$'s be obtained from $t$-tests or that the individual cases correspond to genes. Each of the $N$ cases might involve a separate linear regression for example, with the $i$th case yielding $p$-value $p_i$ for some parameter of interest, and $z_i = \Phi^{-1}(p_i)$.

Section 2 reviews Fdr theory with an emphasis on the *local false discovery rate*, defined in a Bayesian sense as

(1.2) $$\text{fdr}(z_i) = \text{Prob}\{\text{gene } i \text{ is null}|z_i = z\}.$$

An estimate of Fdr $(z)$ for the prostate data is shown by the solid curve in Figure 2, constructed as in Section 2, where it is suggested that a reasonable threshold for reporting likely nonnull genes is $\text{fdr}(z_i) \leq 0.2$. 51 of the 6033 genes have fdr $\leq 0.2$, 25 on the left and 26 on the right. A list of these genes could be reported to the investigators with assurances that it contains less than 20% null cases. Here fdr methods are being used to control *size*, or Type I errors.

The solid bars in Figures 1 and 2 are estimates of the *nonnull histogram*, what we would see if we had $z$ values only for the nonnull genes, constructed as in Section 3. Combined with the fdr curve, the nonnull histogram helps assess *power*, the ability of the data to identify nonnull genes. Figure 2 suggests low power for the prostate data: most of the nonnull cases have large values of $\text{fdr}(z_i)$, and cannot be reported on a list of interesting genes without also reporting a large percentage of null cases. Section 3 constructs some simple power diagnostics based on fdr considerations.

Following [3], most of the Fdr literature has focussed on *tail area false discovery rates*,

(1.3) $$\text{Fdr}(z_i) = \text{Prob}\{\text{gene } i \text{ is null}|z_i \leq z\}$$

(or $\text{Prob}\{\text{null}|z_i \geq z\}$ depending on the sign of $z$). Section 2 discusses the relationship of fdr to Fdr, with relative efficiency calculations presented in Section 5. Local fdr's fit in better with empirical Bayes development, and are featured here, but most of the ideas apply just as well to tail area Fdr's.

The discussion in Sections 2 and 3 assumes that the appropriate null distribution is known to the statistician, perhaps being the theoretical $N(0,1)$ null suggested by (1.1), or its permutation-based equivalent [also nearly $N(0,1)$ for both data sets in Figure 1]. This is tenable for the prostate data, but not for the HIV data. Sections 4 and 5 consider the more difficult and common large-scale testing situation where there is evidence against the theoretical null. Efron [8, 9, 10] discusses estimating an *empirical null* in situations like that for the HIV data where the central histogram does not match $N(0,1)$. Some methodology for constructing empirical nulls is described in Section 4, and its efficiency investigated in Section 5. [It gives empirical null $N(-0.011, 0.75^2)$ for the HIV data, as shown in Figure 1.]

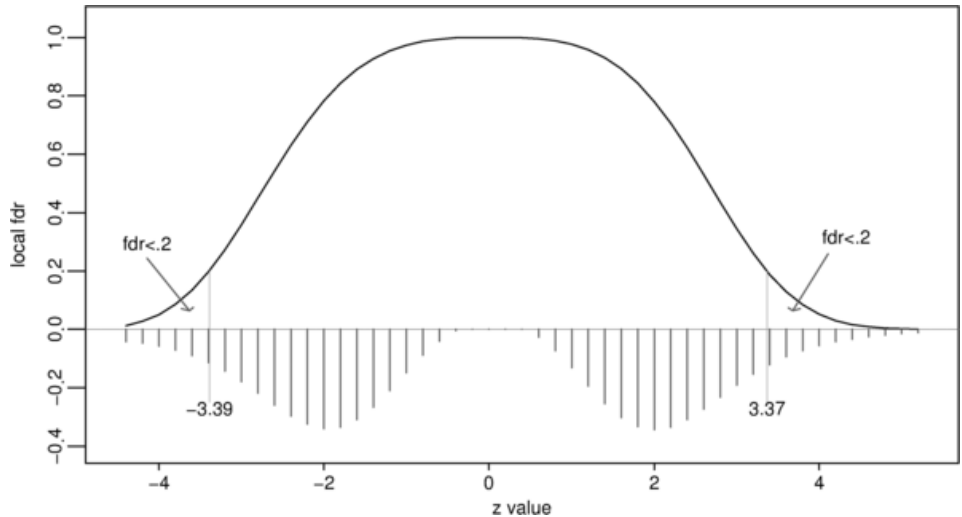

FIG. 2. *Local false discovery rate* fdr($z$), (1.2) *for prostate data, solid curve.* 51 *genes,* 25 *on the left and* 26 *on the right, have* fdr($z_i$) $\leq 0.2$, *a reasonable threshold for reporting nonnull cases. Solid bars show estimated nonnull histogram (plotted negatively, divided by* 50), *constructed as in Section* 3. *Most of the nonnull cases will not be reported.*

Three pairs of related ideas are considered here:

- Size and power calculations for large-scale simultaneous testing.
- Local and tail-area false discovery rates.
- Theoretical and empirical null hypotheses.

All combinations are possible, a power analysis using local fdr with a theoretical null distribution for instance, but only a few are illustrated in the text.

A substantial microarray statistics literature has developed in the past few years, much of it focused on the control of frequentist Type I errors; see, for example, [7], and the review article by Dudoit, Shaffer and Boldruck [6]. Bayes and empirical Bayes methods have also been advocated, as in [18, 19] and [27], while Benjamini and Hochberg's Fdr theory is increasingly influential; see [15] and [33]. Lee et al. [22] and Kerr, Martin and Churchill [20] discuss large-scale inference from ANOVA viewpoints. Local fdr methods, which this article argues can play a useful role, were introduced in [14]; several references are listed at the end of Section 2. The paper ends with a brief summary in Section 6.

**2. False discovery rates.** Local false discovery rates [13, 14], are a variant of [3] "tail area" false discovery rates. This section relates the two ideas, reviews a few basic properties, and presents some general guidelines for interpreting fdr's.

Suppose we have $N$ null hypotheses to consider simultaneously, each with its own test statistic,

$$
\begin{aligned}
\text{Null hypothesis:} \quad & H_1, H_2, \ldots, H_i, \ldots, H_N, \\
\text{Test statistic:} \quad & z_1, z_2, \ldots, z_i, \ldots, z_N;
\end{aligned}
\tag{2.1}
$$

$N$ must be large for local fdr calculations, at least in the hundreds, but the $z_i$ need not be independent. (At least not for getting consistent Fdr estimates, though correlations can decrease the accuracy of such estimates, as detailed in Section 5.) A simple Bayesian "two-class" model, [14, 22, 26], underlies the theory: we assume that the $N$ cases are divided into two classes, null or nonnull, occurring with prior probabilities $p_0$ or $p_1 = 1 - p_0$, and with the density of the test statistic $z$ depending upon its class,

$$
\begin{aligned}
p_0 &= \Pr\{\text{null}\}, & f_0(z) \text{ density if null,} \\
p_1 &= \Pr\{\text{nonnull}\}, & f_1(z) \text{ density if nonnull.}
\end{aligned}
\tag{2.2}
$$

It is natural to take $f_0(z)$ to be a standard $N(0,1)$ density in context (1.1), the *theoretical null*. Here and in Section 3 we assume that $f_0(z)$ is known to the statistician, deferring until Section 4 its estimation in situations like that for the HIV data where the theoretical null is not believable. Fdr theory does not require specification of $f_1(z)$, which is only assumed to be longer-tailed than $f_0(z)$, with the nonnull $z_i$'s tending to occur farther away from 0. Proportion $p_0$, the Bayes a priori probability of a gene being null, is also supposed known here, its estimation being discussed in Sections 4 and 5. Practical applications of large-scale testing usually assume $p_0$ large, say

$$p_0 \geq 0.9, \tag{2.3}$$

the goal being to identify a relatively small set of interesting nonnull cases. Under assumption (2.3), $p_0$ has little practical effect on the usual false discovery rate calculations, that is, on the control of Type I errors, but it will become more crucial for the power diagnostics of Section 3.

Define the *null subdensity*

$$f_0^+(z) = p_0 f_0(z) \tag{2.4}$$

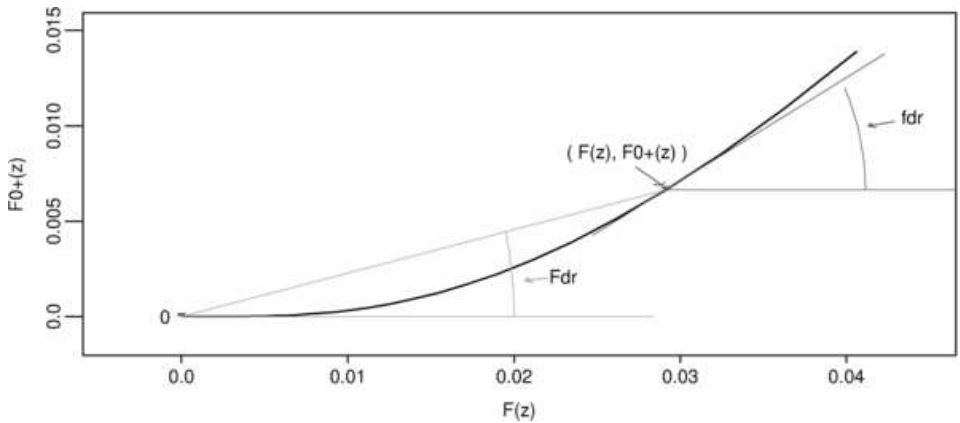

FIG. 3. *Geometrical relationship of Fdr to fdr; heavy curve plots $F_0^+(z)$ versus $F(z)$; $\mathrm{fdr}(z)$ is slope of tangent, $\mathrm{Fdr}(z)$ slope of secant.*

and the *mixture density*

$$f(z) = p_0 f_0(z) + p_1 f_1(z). \tag{2.5}$$

The Bayes posterior probability that a case is null given $z$, by definition the *local false discovery rate*, is

$$\mathrm{fdr}(z) \equiv \Pr\{\mathrm{null}|z\} = p_0 f_0(z)/f(z) = f_0^+(z)/f(z). \tag{2.6}$$

The Benjamini–Hochberg false discovery rate theory relies on tail areas rather than densities. Letting $F_0(z)$ and $F_1(z)$ be the c.d.f.'s corresponding to $f_0(z)$ and $f_1(z)$ in (2.2), define $F_0^+(z) = p_0 F_0(z)$ and $F(z) = p_0 F_0(z) + p_1 F_1(z)$. Then the posterior probability of a case being null given that its $z$-value "$Z$" is less than some value $z$ is

$$\mathrm{Fdr}(z) \equiv \Pr\{\mathrm{null}|Z \leq z\} = F_0^+(z)/F(z). \tag{2.7}$$

(It is notationally convenient to consider events $Z \leq z$ but we could just as well consider tail areas to the right, two-tailed events, etc.) Figure 3 illustrates the geometrical relationship between Fdr and fdr.

Analytically, Fdr is a conditional expectation of fdr [13],

$$\begin{aligned}
\mathrm{Fdr}(z) &= \int_{-\infty}^{z} \mathrm{fdr}(Z) f(Z) \, dZ \Big/ \int_{-\infty}^{z} f(Z) \, dZ \\
&= E_f\{\mathrm{fdr}(Z)|Z \leq z\},
\end{aligned} \tag{2.8}$$

"$E_f$" indicating expectation with respect to $f(z)$ [13]. That is, $\mathrm{Fdr}(z)$ is the average of $\mathrm{fdr}(Z)$ for $Z \leq z$; $\mathrm{Fdr}(z)$ will be less than $\mathrm{fdr}(z)$ in the usual situation where $\mathrm{fdr}(z)$ decreases as $|z|$ gets large. For example $\mathrm{fdr}(-3.39) =$ 0.20 in Figure 2 while $\mathrm{Fdr}(-3.39) = 0.105$. If the c.d.f.'s $F_0(z)$ and $F_1(z)$ are Lehmann alternatives,

$$F_1(z) = F_0(z)^\alpha \qquad [\alpha < 1], \tag{2.9}$$

it is straightforward to show that

$$\log\left\{\frac{\mathrm{fdr}(z)}{1 - \mathrm{fdr}(z)}\right\} = \log\left\{\frac{\mathrm{Fdr}(z)}{1 - \mathrm{Fdr}(z)}\right\} + \log\left(\frac{1}{\alpha}\right), \tag{2.10}$$

giving

$$\mathrm{fdr}(z) \doteq \mathrm{Fdr}(z)/\alpha \tag{2.11}$$

for small values of Fdr. The prostate data of Figure 1 has $\alpha$ roughly $1/2$ in each tail.

Benjamini and Hochberg's [3] Fdr control rule depends on an estimated version of (2.7) where $F(z)$ is replaced by the empirical c.d.f. "$\bar{F}$" of the $z$ values,

$$\widehat{\mathrm{Fdr}}(z) = p_0 F_0(z)/\bar{F}_0(z) \qquad [\bar{F}(z) = \#\{z_i \leq z\}/N]. \tag{2.12}$$

Storey [32] and Efron and Tibshirani [13] discuss the connection of the frequentist Fdr control procedure with Bayesian form (2.7). $\mathrm{Fdr}(z)$ corresponds to Storey's "$q$-value," the tail-area false discovery rate attained at a given observed value $z_i = z$. $\widehat{\mathrm{Fdr}}(z)$ is biased upward as an estimate of $\mathrm{Fdr}(z)$; see Section 4 of [13].

The estimated fdr curve in Figure 2 is

$$\widehat{\mathrm{fdr}}(z) = p_0 f_0(z)/\widehat{f}(z), \tag{2.13}$$

where $f_0(z)$ is the standard normal density $\varphi(z) = \exp\{-z^2/2\}/\sqrt{2x}$, $p_0 = 0.932$ is the value derived in Section 4, and $\widehat{f}(z)$ is a maximum likelihood estimate (MLE) of the mixture density $f(z)$, (2.5), within the seven-parameter exponential family described in Section 4. This type of flexible parametric modeling takes advantage of the fact that $f(z)$, as a mixture of null and nonnull $z$ values, tends to be quite smooth; see Section 6 of [9]. Lindsey's method, Lindsey [24, 25] described in Section 4, finds $\widehat{f}(z)$ using standard Poisson GLM software. The theory and simulations of Section 5 show only a moderate cost in estimation variability for $\widehat{\mathrm{fdr}}$ compared to $\widehat{\mathrm{Fdr}}$.

A variety of other local fdr estimation methods have been suggested: using more specific parametric models such as normal mixtures, see [1, 28, 30] or [17]; isotonic regression [4]; local smoothing [2]; and hierarchical Bayes analyses [5, 23]. All seem to perform reasonably well. The Poisson GLM methodology of this paper has the advantage of easy implementation with familiar software, and permits a closed-form error analysis as shown in Section 5.

Classical frequentist hypothesis testing methods rely on tail areas by necessity. Large-scale testing situations allow us to do local calculations, which are more natural from a Bayesian point of view. For example, the 25 prostate data genes having $z_i \leq -3.39$ have $q$-value $\widehat{\text{Fdr}}(-3.39) = 0.105$; they have average $\widehat{\text{fdr}}(z_i)$ of about 0.105 [as in (2.8)], but varying from 0.20 at the boundary point $z_i = -3.39$ down to $\text{fdr}(z_i) = 0.01$ at $z_i = -4.4$. This just says the obvious, that $z_i$'s further from the boundary are less likely to be false discoveries, which is the useful message conveyed by $\widehat{\text{fdr}}(z)$. The power diagnostics of Section 3 rely on the local Bayesian interpretation (2.6).

The literature has not reached consensus on a standard choice of $q$ for Benjamini–Hochberg testing, the equivalent of 0.05 for single tests, but Bayesian calculations offer some insight. The cutoff threshold $\text{fdr} \leq 0.20$ used in Figure 2 yields posterior odds ratio

$$\Pr\{\text{nonnull}|z\}/\Pr\{\text{null}|z\} = (1 - \text{fdr}(z))/\text{fdr}(z) \tag{2.14}$$
$$= p_1 f_1(z)/p_0 f_0(z) \geq 0.8/0.2 = 4.$$

If we assume prior odds ratio $p_1/p_0 \leq 0.1/0.9$ as in (2.3), then (2.12) corresponds to the Bayes factor

$$f_1(z)/f_0(z) \geq 36 \tag{2.15}$$

in favor of nonnull.

This threshold requires a much stronger level of evidence against the null hypothesis then in standard one-at-a-time testing, where the critical threshold lies somewhere near 3 [11]. We might justify (2.15) as being conservative in guarding against multiple testing fallacies. More pragmatically, increasing the fdr threshold much above 0.20 can deliver unacceptably high proportions of false discoveries to the investigators. The 0.20 threshold, used in the remainder of the paper, corresponds to $q$-values between 0.05 and 0.15 for reasonable choices of $\alpha$ in (2.11); such $q$-value thresholds can be interpreted as reflecting a conservative Bayes factor for Fdr interpretation.

*Any* choice of threshold is liable to leave investigators complaining that the statisticians' list of nonnull cases omits some of their a priori favorites. Conveying the full list of values $\text{fdr}(z_i)$, not just those for cases judged nonnull, allows investigators to employ their own prior opinions on interpreting significance. This is particularly important for low-powered situations like the prostate study, where luck plays a big role in any one case's results, but it is the counsel of perfection, and most investigators will require some sort of reduced list.

The basic false discovery rate idea is admirably simple, and in fact does not depend on the literal validity of the two-class model (2.2). Consider the 28 genes in the prostate example that have $z_i \geq 3.3$; the expected number of null genes having $z_i \geq 3.3$ is 2.71 $[= 6033 \cdot 0.932(1 - \Phi(3.3))]$, so

$$\widehat{\text{Fdr}} = 2.71/28 = 0.097. \tag{2.16}$$

The Fdr interpretation is that about one tenth of the 28 genes can be expected to be null, the other nine tenths being genuine nonnull discoveries.

This interpretation does not require independence, nor even all of the null genes to have the same density $f_0(z)$, only that their *average* density behaves like $f_0$. Since the denominator 28 is observed, the nonnull density $f_1(z)$ plays no role. Exchangeability of the 28 cases is the only real assumption, coming into play when we report that each of the 28 genes has the same one tenth probability of being null. The local fdr has an advantage here, since the equivalent exchangeability assumption is made only for genes having the same observed $z$ values. These ideas are examined in Section 4 of [13].

**3. Power diagnostics.** The microarray statistics literature has focussed on controlling Type I error, the false rejection of genuinely null cases. Dudoit, van der Laan and Pollard [7] provides a good review. Local fdr methods can also help assess power, the probability of rejecting genuinely nonnull cases. This section discusses power diagnostics based on $\widehat{\text{fdr}}(z)$, showing, for example, why the prostate study might easily fail to identify important genes. The emphasis here is on diagnostic statistics that are dependable and simple to calculate.

The nonnull density $f_1(z)$ in the two-class model (2.2), unimportant for the "size" calculations of Section 2, plays a central role in assessing power. From (2.5) and (2.6) we obtain

$$p_1 = \int_{-\infty}^{\infty} [1 - \text{fdr}(z)] f(z) \, dz = 1 - p_0 \tag{3.1}$$

and

$$f_1(z) = (1 - \text{fdr}(z)) f(z)/p_1. \tag{3.2}$$

An estimate of $f(z)$ gives $\widehat{\text{fdr}}(z)$ as in (2.13), and then the estimated nonnull density

$$\widehat{f}_1(z) = [1 - \widehat{\text{fdr}}(z)] \widehat{f}(z) \Big/ \int_{-\infty}^{\infty} [1 - \widehat{\text{fdr}}(z')] \widehat{f}(z') \, dz'. \tag{3.3}$$

Power diagnostics are obtained from the comparison of $\widehat{f}_1(z)$ with $\widehat{\text{fdr}}(z)$. The expectation of $\widehat{\text{fdr}}$ under $\widehat{f}_1$, say "$\widehat{E\text{fdr}}_1$," provides a particularly simple diagnostic statistic,

$$\widehat{E\text{fdr}} = \int_{-\infty}^{\infty} \widehat{\text{fdr}}(z) [1 - \widehat{\text{fdr}}(z)] \widehat{f}(z) \, dz \Big/ \int_{-\infty}^{\infty} [1 - \widehat{\text{fdr}}(z)] \widehat{f}(z) \, dz. \tag{3.4}$$

A small value of $\widehat{E\text{fdr}}_1$, perhaps $\widehat{E\text{fdr}}_1 \doteq 0.20$, suggests good power, with a typical nonnull gene likely to show up on a list of interesting candidates for further study. Neither of the examples in Figure 1 demonstrates good power; $\widehat{E\text{fdr}}_1 = 0.68$ for the prostate data and 0.47 for the HIV data (the latter based on the empirical null $\widehat{\text{fdr}}$ estimate of Section 4).

The *nonnull counts* pictured in Figures 1 and 2 allow a more intuitive interpretation of formula (3.4). Suppose that the $N$ $z$-values have been placed into $K$ bins of equal width $\Delta$, with

(3.5)
$$x_k = \text{centerpoint of } k\text{th bin} \quad \text{for } k = 1, 2, \ldots, K,$$
$$y_k = \#\{z_i \text{ in } k\text{th bin}\}.$$

Since

(3.6)  $$\text{Prob}\{\text{gene } i \text{ nonnull}|z_i = z\} = 1 - \text{fdr}(z),$$

an approximately unbiased estimate of the nonnull counts in bin $k$ is

(3.7)  $$\widehat{y}_{1k} = [1 - \widehat{\text{fdr}}(x_k)] \cdot y_k.$$

The solid bars in Figures 1 and 2 follow definition (3.7), except with $y_k$ replaced by a smoothed estimate proportioned to $\widehat{f}(x_k)$. Looking at Figure 2, an obvious estimate of the nonnull expectation $E\text{fdr}_1$ is

(3.8)  $$\widehat{E\text{fdr}}_1 = \frac{\sum_k \widehat{\text{fdr}}(x_k) \cdot \widehat{y}_{1k}}{\sum_k \widehat{y}_{1k}} \doteq \frac{\sum_k \widehat{\text{fdr}}(x_k)[1 - \widehat{\text{fdr}}(x_k)]\widehat{f}(x_k)}{\sum_k [1 - \widehat{\text{fdr}}(x_k)]\widehat{f}(x_k)},$$

which amounts to evaluating the integrals in (3.4) via the trapezoid rule.

Table 1 reports on a simulation study of $\widehat{E\text{fdr}}_1$. The study took

(3.9) $$z_i \stackrel{\text{ind}}{\sim} N(\mu_i, 1) \quad \text{with } \begin{cases} \mu_i = 0, & \text{probability } 0.90, \\ \mu_i \sim N(3, 1), & \text{probability } 0.10, \end{cases}$$

for $i = 1, 2, \ldots, N = 1500$. [More precisely, $\mu_i = 3 + \Phi^{-1}((i - 0.5)/150)$, $i = 1, 2, \ldots, 150$, for the nonnull cases.] The "theoretical null" columns assume $f_0 = N(0, 1)$, while "empirical null" estimates $f_0$ by the central matching method of Section 4. Both methods estimated $\widehat{p}_1 = 1 - \widehat{p}_0$ by central matches. The true value of $E\text{fdr}_1$ in situation (3.9) is 0.32. The estimates $\widehat{E\text{fdr}}_1$ are seen to be reasonably stable and roughly accurate. Section 5 discusses the downward bias in $\widehat{p}_1$.

Going further, we can examine the entire distribution of $\widehat{\text{fdr}}$ under $\widehat{f}_1$ rather than just its expectation. The nonnull c.d.f. of $\widehat{\text{fdr}}$ is estimated by

(3.10)  $$\widehat{G}_1(t) = \sum_{k:\widehat{\text{fdr}}(x_k) \leq t} \widehat{y}_{1k} / \sum_k \widehat{y}_{1k}$$

TABLE 1
*Means, standard deviations and coefficients of variation of $\widehat{E\mathrm{fdr}}_1$ (3.5); 100 trials of situation (3.9), $N = 1500$. True value $E\mathrm{fdr}_1 = 0.32$, $p_1 = 0.10$*

|         | Theoretical null |       | Empirical null |       |
|---------|:---:|:---:|:---:|:---:|
|         | $\widehat{E\mathrm{fdr}}_1$ | $\widehat{p}_1$ | $\widehat{E\mathrm{fdr}}_1$ | $\widehat{p}_1$ |
| Mean    | **0.285** | 0.085 | **0.232** | 0.076 |
| Stdev   | 0.060 | 0.015 | 0.040 | 0.011 |
| Coefvar | 0.21  | 0.18  | 0.17  | 0.14  |

for $0 \le t \le 1$. Figure 3 shows $\widehat{G}_1(t)$ for the prostate study, for the HIV study [taking $f_0 = N(-0.11, 0.75^2)$, $p_0 = 0.93$, as in Figure 1], and for the first of the 100 simulations from model (3.9). The simulation curve suggests good power characteristics, with 64% of the nonnull genes having fdr less than 0.2. At the opposite extreme, only 11% of nonnull genes in the prostate study have fdr less than 0.2.

Graphs such as Figure 4 help answer the researchers' painful question "why are not the genes we expected on your list of nonnull outcomes?" For the prostate data, *most* of the nonnull genes will not turn up on a list of low fdr cases. The R program *locfdr*, discussed in Section 4, returns $\widehat{E\mathrm{fdr}}_1$ and a graph of $\widehat{G}_1(t)$.

Traditional sample size calculations employ preliminary data to predict how large an experiment might be required for effective power. Here we might ask, for instance, if doubling the number of subjects in the prostate study would substantially improve its detection rate.

To answer this question, denote the mean and variance of $z_i$ by $\mu_i$ and $\sigma_i^2$,

$$(3.11) \qquad z_i \sim (\mu_i, \sigma_i^2).$$

We imagine that $c$ independent replicates of $z_i$ are available for each gene (doubling the experiment corresponding to $c = 2$), from which a combined test statistic $\widetilde{z}_i$ is formed,

$$(3.12) \qquad \widetilde{z}_i = \sum_{j=1}^{c} z_{ij}/\sqrt{c} \sim (\sqrt{c}\mu_i, \sigma_i^2).$$

This definition maintains the mean and variance of null cases, $\widetilde{z}_i \sim (0, \sigma_i^2)$, while moving the nonnull means $\widetilde{\mu}_i = \sqrt{c}\mu_i$ away from zero by the factor $\sqrt{c}$.

Consider a subset of $m$ nonnull genes, say $\mathcal{S}$, and define

$$(3.13) \quad \bar{\mu} = \sum_{\mathcal{S}} \mu_i/m, \qquad \Delta^2 = \sum_{\mathcal{S}} (\mu_i - \bar{\mu})^2/m \quad \text{and} \quad \bar{\sigma}^2 = \sum_{\mathcal{S}} \sigma_i^2/m.$$

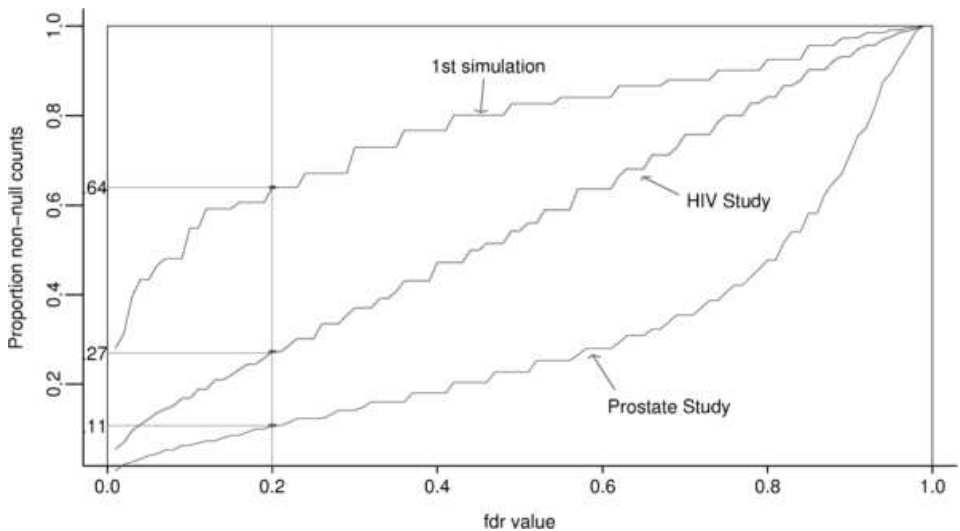

Fig. 4. *Estimated nonnull c.d.f. of fdr, (3.10); prostate study, HIV study, and first of 100 simulations, model (3.9). The simulation curve suggests substantial power, with 64% of the nonnull cases having $\widehat{\text{fdr}}$ less than 0.2. $\widehat{E\text{fdr}}_1$ values: 0.23 (simulation), 0.45 (HIV), 0.68 (prostate).*

A randomly selected $z_i$ value "$Z$" from $\mathcal{S}$ has mean and variance

$$(3.14) \qquad Z \sim (\bar{\mu}, \Delta^2 + \bar{\sigma}^2),$$

while $\widetilde{Z}$, the corresponding randomly selected $\widetilde{z}_i$ value, has

$$(3.15) \qquad \widetilde{Z} \sim (\sqrt{c}\bar{\mu}, c\Delta^2 + \bar{\sigma}^2),$$

so

$$(3.16) \quad \widetilde{Z} = \sqrt{c}\,\bar{\mu} + d(Z - \bar{\mu}) \qquad [d^2 = c - (c-1)\bar{\sigma}^2/(\Delta^2 + \bar{\sigma}^2)]$$

gives $\widetilde{Z}$ the correct mean and variance.

Let $\mathcal{S}$ be the set of nonnull genes having $z_i > 0$. We can estimate $\bar{\mu}$ and $d$ from the corresponding nonnull counts $\widehat{y}_{1k}$ (the bars to the right of $z = 0$ in Figure 2), and then use (3.16) to move those counts out from location $x_k$ to $\widetilde{x}_k = \sqrt{c}\{\bar{\mu} + d^{1/2}(x_k - \bar{\mu})\}$. The null counts $\widehat{y}_{0k} = y_k - \widehat{y}_{1k}$ do *not* change location when the sample size increases. We can do the same calculations for the nonnull counts having $z_i < 0$. Together, these provide an estimate of what the entire $z$-value histogram would look like if the sample size were increased by the factor $c$, from which we can recalculate $\widehat{E\text{fdr}}_1$ or any other diagnostic statistic.

Table 2 shows $\widehat{E\text{fdr}}_1$ estimates for hypothetical expansions of the prostate and HIV studies. Doubling the HIV study, to 8 instead of 4 subjects in each HIV category, reduces $\widehat{E\text{fdr}}_1$ from 0.45 to 0.23, while doubling the prostate study gives less dramatic improvement. Table 2 is based on a cruder version of (3.16) that takes $d = 1$,

$$(3.17) \qquad \widetilde{Z} = \sqrt{c}Z,$$

in other words, simply moving the nonnull counts $\widehat{y}_{1k}$ from $x_k$ to $\sqrt{c}x_k$. Using (3.17) tends to underestimate the reduction in $\widehat{E\text{fdr}}_1$, but did not make much difference in this case.

TABLE 2
*Hypothetical values of $\widehat{E\text{fdr}}_1$ for versions of the prostate and HIV studies expanded by factor c; based on transformation (3.17) for the nonnull counts, as calculated by R program* locfdr

| $c$ | 1 | 1.5 | 2 | 2.5 | 3 |
|---|---|---|---|---|---|
| Prostate | 0.68 | 0.54 | 0.44 | 0.38 | 0.34 |
| HIV | 0.45 | 0.31 | 0.23 | 0.18 | 0.14 |

The R program *locfdr* does these calculations. They have a speculative nature, but no more so than traditional power projections. Like all of the diagnostics in this section, they require no mathematical assumptions beyond the original two-class model (2.2).

**4. Empirical null estimation.** The null density $f_0(z)$ in (2.2) is crucial to false discovery rate calculations, or for that matter to any hypotheses testing method. We assumed $f_0 \sim N(0,1)$, the *theoretical null*, for the prostate data. This seems natural in situation (1.1), being almost certainly what would be done if there were only a single gene's data to analyze. Large scale simultaneous testing, however, raises the possibility of detecting deficiencies in the theoretical null, as with the HIV data in Figure 1 where the $z$-value histogram is noticeably too narrow around its central peak. This section concerns data-based estimates of $f_0(z)$, for example, the *empirical null distribution* $\widehat{f}_0 \sim N(-0.11, 0.75^2)$ for the HIV data, shown in Figure 1.

Efron [8, 10] lists four reasons why $f_0$ might differ from the theoretical null:

(1) *Failed assumptions.* Let $Y$ be the $N$ by $n$ matrix of expression levels, $N$ genes and $n$ microarrays in our two studies,

$$(4.1) \qquad y_{ij} = \text{expression level for } i\text{th gene and } j\text{th array.}$$

The HIV study has $N = 7680$ genes and $n = 8$ microarrays, 4 each from HIV positive and HIV negative subjects. Each gene yields a two-sample $t$ statistic $t_i$ comparing positive versus negative subjects, with $z$-value

$$(4.2) \qquad z_i = \Phi^{-1}(F_6(t_i)),$$

$F_6$ the c.d.f. of a $t$ distribution having 6 degrees of freedom.

The theoretical null distribution $f_0 \sim N(0,1)$ is justified for (4.2) if the $y_{ij}$'s are normal, or by asymptotic theory as $n$ goes to infinity, neither argument applying to the HIV data. We can avoid these assumptions by computing a *permutation null*, the marginal distribution of the $z_i$'s obtained by permuting the columns of $Y$. This gave $\widehat{f}_0 \dot{\sim} N(0, 0.99^2)$ for the HIV data, failing to explain the narrow central peak in Figure 1.

(2) *Unobserved covariates.* The HIV study is observational: subjects were observed, not assigned to be HIV positive or negative, and similarly for the prostate study. Section 4 of [8] discusses how unobserved covariates in observational studies are likely to widen $f_0(z)$, and how this effect is not detectable by permutation analysis. A microarray example is presented in which the $z$-value histogram has central dispersion more than half again as wide as the theoretical null. Since the HIV histogram is too narrow at its center rather than too wide, unobserved covariates are not the problem here.

(3) *Correlation across arrays.* The theoretical null as applied to (4.2) or (1.1) assumes independence across the columns of $Y$, that is, among $y_{i1}, y_{i2}, \ldots, y_{in}$ for gene $i$. Experimental difficulties can undercut independence in microarray studies, while being undetectable in the permutation null distribution. The HIV data was analyzed with the HIV negative subjects as the first four columns of $Y$ and the positives as the last four columns. A principal components analysis suggested a strong pattern of correlation *across* columns, with arrays $(1,3,5,7)$ positively correlated, and likewise arrays $(2,4,6,8)$. This pattern would narrow the null distribution in situation (4.2).

(4) *Correlation between genes.* A striking advantage of the two-group model and its false discovery rate analysis in Section 2 is that it does *not* require independence between genes. Estimates such as $\widehat{\mathrm{fdr}}(z) = p_0 f_0(z)/\widehat{f}(z)$ only require consistency for $\widehat{f}(z)$ (but do not achieve the full efficiency attainable from knowledge of the gene-wise correlation structure).

Efron [10] emphasizes a caveat to this argument: even if the theoretical null is individually appropriate for each null gene, correlations between genes can make the effective null distribution $f_0(z)$ substantially narrower or wider than $N(0,1)$. There it is shown that the amount of correlation in the HIV data could easily explain a $N(-0.11, 0.75^2)$ null distribution. (By contrast, the prostate data exhibits quite small gene-wise correlations.) A permutation null distribution will not reveal correlation effects.

Empirically estimating the null distribution avoids all four difficulties, and any others that may distort $f_0$. There is a price to pay, though, in terms of accuracy: using the empirical null substantially increases the variability of estimated false discovery rates, as shown in Section 5. This price is unavoidable in situations like the HIV study where there is clear evidence against the theoretical null; the null distribution provides the crucial numerator in false discovery rate estimates like (2.16), where using an inappropriate null undercuts inferential validity. (Using the theoretical null on the HIV data eliminates all but 20 of the 151 genes having empirical null $\widehat{\text{fdr}}$ estimate less than 0.20, including all of those with $z_i < 0$.)

The basic empirical null idea is simple: we assume $f_0(z)$ is normal but not necessarily with mean 0 and variance 1, say

(4.3) $$f_0(z) \sim N(\delta_0, \sigma_0^2),$$

and then estimate $\delta_0, \sigma_0$, as well as the null proportion $p_0$ in (2.2), from the histogram data near $z = 0$. Two different methods for estimating $(\delta_0, \sigma_0, p_0)$ will be described, and their accuracies analyzed in Section 5. Both methods are implemented in algorithm *locfdr*, available through the Comprehensive R Archive Network, http://cran.r-project.org; *locfdr* produced all of the numerical examples in this paper.

"Central matching," the first of our two estimation methods for $(\delta_0, \sigma_0, p_0)$, operates by quadratically approximating $\log \widehat{f}(z)$ around $z = 0$. To begin, the *locfdr* algorithm estimates $f(z)$, (2.5), by maximum likelihood fitting to the histogram counts $y_k$ for the $z$ values, (3.5), within a parametric exponential family. Figure 2 used the seven-parameter family

(4.4) $$f_{\boldsymbol{\beta}}(\mathbf{z}) = c_{\boldsymbol{\beta}} \exp\left\{\sum_{j=1}^{7} \beta_j z^j\right\},$$

$c_{\boldsymbol{\beta}}$ the constant making $f_{\boldsymbol{\beta}}$ integrate to 1.

Figure 5 illustrates central matching estimation of $(\delta_0, \sigma_0, p_0)$ for the HIV data based on the methodology in [8]. The heavy curve is $\log \widehat{f}(z)$, fit by maximum likelihood [using a natural spline basis with 7 degrees of freedom, rather than the polynomials of (4.4), another option in *locfdr*, though (4.4) gives nearly identical results in this case]. A quadratic curve $\widehat{f}_0^+(z)$ has been fit to $\log \widehat{f}(z)$ around $z = 0$,

(4.5) $$\log(\widehat{f}_0^+(z)) = \widehat{\beta}_0 + \widehat{\beta}_1 z + \widehat{\beta}_2 z^2.$$

Assuming $f_0(z) \sim N(\delta_0, \sigma_0^2)$, the log of the null subdensity (2.4) is

(4.6) $$\log(f_0^+(z)) = \log p_0 - \frac{1}{2}\left\{\frac{\delta_0^2}{\sigma_0^2} + \log(2\pi\sigma_0^2)\right\} + \frac{\delta_0}{\sigma_0^2}z - \frac{1}{2\sigma_0^2}z^2.$$

Estimates $(\widehat{\beta}_0, \widehat{\beta}_1, \widehat{\beta}_2)$ from (8.5) translate to estimates $(\widehat{\delta}_0, \widehat{\sigma}_0, \widehat{p}_0)$ in (4.6), for example, $\widehat{\sigma}_0 = (2\widehat{\beta}_2)^{-1/2}$. For the HIV data this gave

(4.7) $$\widehat{\delta}_0 = -0.107, \quad \widehat{\sigma}_0 = 0.753 \quad \text{and} \quad \widehat{p}_0 = 0.931.$$

The logic here is straightforward: we make the "zero assumption" that the central peak of the $z$-value histogram consists mainly of null cases, and choose $(\delta_0, \sigma_0, p_0)$ in (4.6) to quadratically approximate the histogram counts near $\delta = 0$. Some form of the zero assumption is required because the two-class model (2.2) is unidentifiable in the absence of strong parametric assumptions on $f_1$.

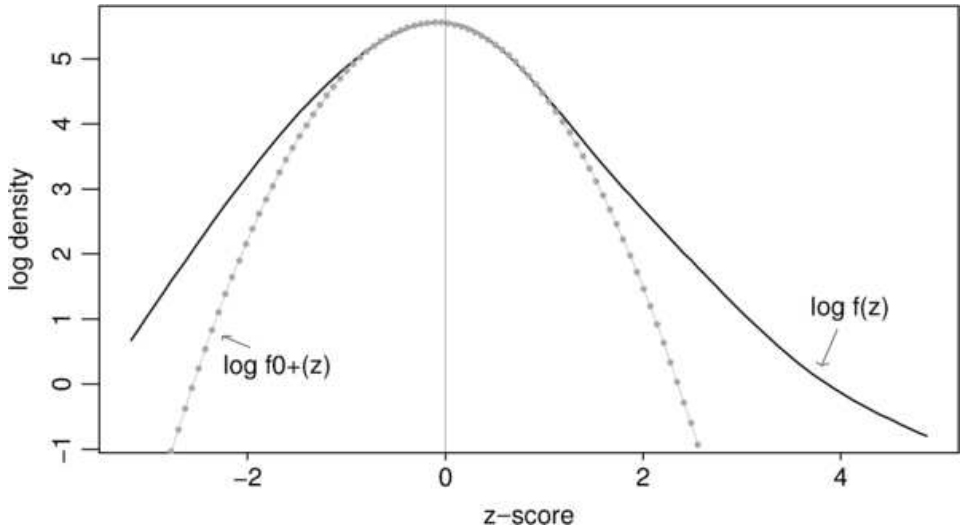

FIG. 5. *Central matching estimation of $p_0$ and $f_0(z) \sim N(\delta_0, \sigma_0^2)$ for the HIV data; heavy curve is log of $\widehat{f}(z)$, estimated mixture density* (2.5); *beaded curve is quadratic fit to* $\log \widehat{f}(z)$ *around* $z = 0$, *estimating* $\log f_0^+(z)$, (2.4). *The three estimated coefficients of quadratic fit give* $(\widehat{\delta}_0, \widehat{\sigma}_0, \widehat{p}_0)$.

A healthy literature exists on estimating $p_0$, as in [21] and [29], all of which relies on the zero assumption [mostly working with $p$-values rather than $z$-values, e.g., $p_i = F_6(t_i)$ in (4.2), where the "zero region" occurs near $p = 1$]. All of this literature relies on the validity of the theoretical null, so in this sense (4.5) and (4.6) is a straightforward extension to situations where the theoretical null is untrustworthy. For the HIV data, using the theoretical null in (4.5) and (4.6), that is, taking $(\widehat{\beta}_1, \widehat{\beta}_2)$ equal $(0, 1/2)$, results in the impossible estimate $\widehat{p}_0 = 1.18$. This will always happen when the $z$-value histogram is narrower than $N(0, 1)$ near $z = 0$.

The zero assumption is more believable if $p_0$, the proportion of null cases, is large. Efron [8] shows that if $p_0$ exceeds 0.90 the fitting method in Figure 5 will be nearly unbiased: although the 10% or less of nonnull cases might in fact contribute some counts near $z = 0$, they cannot substantially affect $\widehat{\delta}_0$ and $\widehat{\sigma}_0$; the $p_0$ estimate *is* affected, being upwardly biased, as seen in Table 1.

"MLE fitting," the second, newer, method of estimating $(\delta_0, \sigma_0, p_0)$, is based on a truncated normal model. We assume that the nonnull density is supported outside some known interval $[-x_0, x_0]$,

(4.8) $$f_1(z) = 0 \quad \text{for } z \in [-x_0, x_0].$$

We need the following definitions:

(4.9) 
$$\mathcal{I}_0 = \{i : z_i \in [-x_0, x_0]\},$$
$$N_0 = \text{number of } z_i \in [-x_0, x_0],$$
$$\mathbf{z}_0 = \{z_i, i \in \mathcal{I}_0\},$$
$$H_0(\delta_0, \sigma_0) = \Phi\left(\frac{x_0 - \delta_0}{\sigma_0}\right) - \Phi\left(\frac{-x_0 - \delta_0}{\sigma_0}\right)$$

and

(4.10) $$\varphi_{\delta_0,\sigma_0}(z) = \frac{1}{\sqrt{2\pi\sigma_0^2}} \exp\left\{-\frac{1}{2}\left(\frac{z - \delta_0}{\sigma_0}\right)^2\right\}.$$

Then

(4.11) $$\theta \equiv p_0 H_0(\delta_0, \sigma_0) = \text{Prob}\{z_i \in [-x_0, x_0]\}$$

under model (2.2).

The likelihood function of the data $(N, \mathbf{z}_0)$ is

(4.12) $$f_{\delta_0,\sigma_0,p_0}(N, \mathbf{z}_0) = [\theta^{N_0}(1-\theta)^{N-N_0}]\left[\prod_{\mathcal{I}_0} \frac{\varphi_{\delta_0,\sigma_0}(z_i)}{H_0(\delta_0, \sigma_0)}\right].$$

This is a product of two exponential family likelihoods, as discussed in Section 5. It is easy to numerically find the MLE estimates $(\widehat{\delta}_0, \widehat{\sigma}_0, \widehat{\theta})$ in (4.12), after which

(4.13) $$\widehat{p}_0 = \widehat{\theta}/H_0(\widehat{\delta}_0, \widehat{\sigma}_0).$$

Table 3 compares the estimates $(\widehat{\delta}_0, \widehat{\sigma}_0, \widehat{p}_0)$ obtained from central matching and MLE fitting for the same 100 simulations of model (3.9) used in Table 1. MLE fitting does somewhat better overall, especially for $\widehat{\delta}_0$, the mean of $\widehat{f}_0$. The results are encouraging, in particular showing that $\sigma_0$ can be estimated within a few per cent. Delta method formulas for the standard deviations are developed in Section 5. These performed well, giving nearly the correct average values, as shown in the table, with small coefficient of variation across the 100 simulations, about 10% for central matching and 3% for MLE fitting. Changing sample size $N = 1500$ by multiple "$c$" changes the standard deviations by about $1/\sqrt{c}$.

The $z_i$'s are independent in model (3.9). This is unrealistic for microarray studies, but as discussed in Section 5, the results may still be applicable to highly correlated situations.

TABLE 3
*Comparison of estimates $(\widehat{\delta}_0, \widehat{\sigma}_0, \widehat{p}_0)$, central matching and MLE fitting; 100 simulations, model (3.9), as in Table 1. MLE fitting took $x_0 = 2$ in (4.8); "formula" standard deviations from delta method calculations, Section 5. True values $(\delta_0, \sigma_0, p_0) = (0, 1, 0.9)$*

|  | Central matching | | | MLE fitting | | |
|---|---|---|---|---|---|---|
|  | mean | stdev | (formula) | mean | stdev | (formula) |
| $\widehat{\delta}_0$ | 0.021 | **0.056** | (0.062) | 0.044 | **0.031** | (0.032) |
| $\widehat{\sigma}_0$ | 1.020 | **0.029** | (0.033) | 1.035 | **0.031** | (0.031) |
| $\widehat{p}_0$ | 0.924 | **0.013** | (0.015) | 0.933 | **0.009** | (0.011) |

The two fitting methods have different virtues and defects. Central matching is attractive from a theoretical point of view, suggesting how we might go beyond normality assumption (4.3), as discussed in Section 7 of [9]. As mentioned before, it gives nearly unbiased estimates of $\widehat{\delta}_0$ and $\widehat{\sigma}_0$ if $p_0$ exceeds 0.9. However, it can be excessively variable, especially in estimating $\delta_0$, and is sensitive to the range of discretization (though not the grid size $\Delta$) in (3.5); reducing the range of $x_k$ in Table 3 from $[-4, 7.4]$ to $[-4, 6.1]$ gave notably worse performance.

MLE fitting generally gives more stable parameter estimates, for reasons suggested by the influence function analysis of Section 5. It does not require discretization of the $z$-values. It does, however, depend strongly on the choice of $x_0$ in (4.8), which was arbitrarily set at $x_0 = 2$ in the simulations. A more adaptive version that began by estimating an appropriate "zero assumption" interval is feasible but more variable. This contrasts with central matching, which automatically adjusts to each situation, including ones where $(\delta_0, \sigma_0)$ is far from $(0, 1)$.

*Locfdr* defaults to MLE fitting for fdr estimation, but also returns the central matching estimates. The two methods gave similar results for the HIV data, $(\widehat{\delta}_0, \widehat{\sigma}_0, \widehat{p}_0) = (-0.117, 0.785, 0.955)$ for MLE fitting, compared to (4.7).

Accurate estimation of $(\delta_0, \sigma_0, p_0)$ is just as important for tail-area Fdr analysis (2.7) as for the local version (2.6). Section 5 computes the accuracy of both $\widehat{\text{Fdr}}(z)$ and $\widehat{\text{fdr}}(z)$. Using an empirical null is expensive in either venue, but the theoretical null can be an unrealistic option, as for the HIV data.

Permutation methods are popular in the microarray literature, but they only address the first of our four listed difficulties for the theoretical null; in practice permutation null distributions are usually close to the theoretical distribution, especially for $t$-test statistics. Permutation and empirical null methods can be used together: if $\widetilde{F}(t)$ is the permutation c.d.f. for the $t$- statistics in the HIV study, then we could set $z_i = \Phi^{-1}(\widetilde{F}(t))$ rather than (4.2) to begin empirical null estimation.

**5. Influence and accuracy.** Two different classifications of false discovery rate methods have been discussed: local versus tail area definitions, Section 2, and theoretical versus empirical null estimates, Section 4. This section derives accuracy formulas for all four combinations, based on closed-form influence functions. In Figure 2, for example, the R algorithm *locfdr* reports $\widehat{\text{fdr}}(3.37) = 0.2 \pm 0.02$ for the theoretical null and local fdr combinations, the standard error 0.02 coming from Theorem 1 below. The influence functions also help explicate differences between central matching and MLE fitting estimation for the empirical null.

For numerical calculations it is convenient to assume that the $N$ $z$-values have been binned as in (3.5): into $K$ bins of width $\Delta$, centerpoint $x_k$, for $k = 1, 2, \ldots, K$, with $y_k$ the count in bin $k$. "Lindsey's method," as discussed in Section 2 of [12], then permits almost fully efficient parametric density estimation within exponential families such as (4.4), using standard Poisson regression software.

*Locfdr* first fits an estimated mixture density $\widehat{f}(z)$, (2.5), to the count vector $\mathbf{y} = (y_1, y_2, \ldots, y_K)'$, by maximum likelihood estimation within a parametric family such as (4.4). For central matching, $\log \widehat{f}_0^+(z)$, (2.4), is fit to $\log \widehat{f}(z)$ as in (4.5)–(4.7); the fitting is by ordinary least squares over a central subset of $K_0$ bins having index set say "$\mathbf{i}_0$." In Figure 5, $\widehat{f}(z)$ was estimated using $K = 41$ bins having centerpoints $-4.0, -3.8, \ldots, 4.0$, while $\log \widehat{f}_0^+(z)$ was fit to $\log \widehat{f}(z)$ from the $K_0 = 6$ central bins $\mathbf{i}_0 = (18, 19, 20, 21, 22, 23)$.

Let $X$ be the $K \times m$ structure matrix used for estimating $\log \widehat{f}(z)$; $X$ has $m = 8$, $k$th row $(1, x_k, x_k^2, \ldots, x_k^7)$ in (4.4). Also let $X_0$ be the $K \times m_0$ matrix used to describe $\log \widehat{f}_0(z)$; $X_0$ has $k$th row $(1, x_k, x_k^2)$, $m_0 = 3$ for the empirical null estimate (4.5), while $X_0$ is the $K \times 1$ matrix $(1, 1, \ldots, 1)'$ for the theoretical null. (Section 7 of [9] considers more ambitious empirical null estimates, e.g., including a cubic term.)

Define submatrices of $X$ and $X_0$,

(5.1) $$\widetilde{X} = X[\mathbf{i}_0, \cdot] \quad \text{and} \quad \widetilde{X}_0 = X_0[\mathbf{i}_0, \cdot],$$

of dimensions $K_0 \times m$ and $K_0 \times m_0$; also

(5.2) $$\widehat{\nu}_k = N\Delta \widehat{f}(x_k), \qquad k = 1, 2, \ldots, K,$$

an estimate of the expected count in bin $k$; and

(5.3) $$\widehat{G} = X' \operatorname{diag}(\widehat{\nu})X, \qquad \widetilde{G}_0 = \widetilde{X}_0' \widetilde{X}_0,$$

where $\operatorname{diag}(\widehat{\nu})$ is a $K \times K$ diagonal matrix having diagonal elements $\widehat{\nu}_k$. Finally, let $\widehat{\boldsymbol{\ell}}$ indicate the $K$-vector with elements $\widehat{\ell}_k = \log \widehat{f}(x_k)$, likewise $\widehat{\boldsymbol{\ell}}_0^+$ for the vector $(\log \widehat{f}_0^+(x_k))$ and $\widehat{\boldsymbol{\ell}\mathbf{fdr}}$ for $(\log \widehat{\text{fdr}}(x_k))$.

By definition the *influence function* of vector $\widehat{\boldsymbol{\ell}\mathbf{fdr}}$ with respect to count vector $\mathbf{y}$ is the $K \times K$ matrix $d\widehat{\boldsymbol{\ell}\mathbf{fdr}}/d\mathbf{y}$ of partial derivatives $\partial \widehat{\ell\mathrm{fdr}}_k/\partial y_\ell$.

LEMMA 1. *The influence function of* $\log \widehat{\mathbf{fdr}}$ *with respect to* $\mathbf{y}$, *when using central matching, is*

$$\text{(5.4)} \qquad \frac{d\widehat{\boldsymbol{\ell}\mathbf{fdr}}}{d\mathbf{y}} = A\widehat{G}^{-1}X',$$

*where*

$$\text{(5.5)} \qquad A = X_0 \widetilde{G}_0^{-1} \widetilde{X}_0' \widetilde{X} - X.$$

PROOF. A small change $d\mathbf{y}$ in the count vector (considered as continuous) produces the change $\mathbf{d}\widehat{\boldsymbol{\ell}}$ in $\widehat{\boldsymbol{\ell}}$,

$$\text{(5.6)} \qquad \mathbf{d}\widehat{\boldsymbol{\ell}} = X\widehat{G}^{-1}X'\,\mathbf{dy}.$$

Similarly if $\widehat{\boldsymbol{\ell}}_0^+ = X_0 \widehat{\gamma}, \widehat{\gamma}$ an $m_0$-vector, is fit by least squares to $\widetilde{\boldsymbol{\ell}} = \widehat{\boldsymbol{\ell}}[\mathbf{i}_0]$, we have

$$\text{(5.7)} \qquad d\widehat{\gamma} = \widetilde{G}_0^{-1} \widetilde{X}_0' \mathbf{d}\widetilde{\boldsymbol{\ell}} \quad \text{and} \quad \mathbf{d}\widehat{\boldsymbol{\ell}}_0^+ = X_0 \widetilde{G}_0^{-1} \widetilde{X}_0' \mathbf{d}\widetilde{\boldsymbol{\ell}}.$$

Both (5.6) and (5.7) are standard regression results. Then (5.6) gives $\mathbf{d}\widetilde{\boldsymbol{\ell}} = \mathbf{d}\widehat{\boldsymbol{\ell}}[\mathbf{i}_0] = \widetilde{X}\widehat{G}^{-1}X'\,\mathbf{dy}$, yielding

$$\mathbf{d}\widehat{\boldsymbol{\ell}}_0^+ = X_0 \widetilde{G}_0^{-1} \widetilde{X}_0' \widetilde{X} \widehat{G}^{-1} X' \,\mathbf{dy}$$

from (5.7). Finally,

$$\text{(5.8)} \qquad \mathbf{d}(\widehat{\boldsymbol{\ell}\mathbf{fdr}}) = \mathbf{d}(\widehat{\boldsymbol{\ell}}_0^+ - \widehat{\boldsymbol{\ell}}) = (X_0 \widehat{G}_0^{-1} \widetilde{X}_0' \widetilde{X} - X)\widehat{G}^{-1} X'\,\mathbf{dy},$$

verifying (5.4). □

THEOREM 1. *In the case where the $z$ values are independent, the delta-method estimate of covariance for the vector of* $\log \widehat{\mathrm{fdr}}(x_k)$ *values, based on central matching, is*

$$\text{(5.9)} \qquad \widehat{\mathrm{cov}}(\widehat{\boldsymbol{\ell}\mathbf{fdr}}) = A\widehat{G}^{-1}A'.$$

PROOF. Under independence, $\mathbf{y}$ has a multinomial distribution with covariance matrix

$$\text{(5.10)} \qquad \mathrm{cov}(\mathbf{y}) = \mathrm{diag}(\boldsymbol{\nu}) - \boldsymbol{\nu}\boldsymbol{\nu}'/N,$$

where $\boldsymbol{\nu} \equiv E\{\mathbf{y}\}$ [$\nu_k \doteq N\Delta f(x_k)$, as in (5.2)]. The delta-method covariance estimate is

$$\text{(5.11)} \qquad \left(\frac{d\widehat{\boldsymbol{\ell}\mathbf{fdr}}}{\mathbf{dy}}\right) \widehat{\mathrm{cov}}(\mathbf{y}) \left(\frac{d\widehat{\boldsymbol{\ell}\mathbf{fdr}}}{\mathbf{dy}}\right)' = (A\widehat{G}^{-1}X)\,\mathrm{diag}(\widehat{\boldsymbol{\nu}})(A\widehat{G}^{-1}X)'$$

$$= A\widehat{G}^{-1}A'.$$

TABLE 4
Boldface: *standard errors of* $\log \widehat{\mathrm{fdr}}(z)$, *("local"), and* $\log \widehat{\mathrm{Fdr}}(z)$, *("tail"),* 250 *replications of model* (3.9), $N = 1500$. *Parentheses: average standard deviation estimate from formula* (5.9); *fdr is the true false discovery rate* (2.6) *based on model* (3.9). *Natural spline basis,* 7 *degrees of freedom used to fit* $f(z)$, *central matching for* $\widehat{f}_0^+(z)$, *empirical null case*

| | | Theoretical null | | | Empirical null | | |
|---|---|---|---|---|---|---|---|
| $z$ | fdr | local | (formula) | tail | local | (formula) | tail |
| 1.5 | 0.88 | **0.05** | (0.05) | **0.05** | **0.04** | (0.04) | **0.10** |
| 2.0 | 0.69 | **0.08** | (0.09) | **0.05** | **0.09** | (0.10) | **0.15** |
| 2.5 | 0.38 | **0.09** | (0.10) | **0.05** | **0.16** | (0.16) | **0.23** |
| 3.0 | 0.12 | **0.08** | (0.10) | **0.06** | **0.25** | (0.25) | **0.32** |
| 3.5 | 0.03 | **0.10** | (0.13) | **0.07** | **0.38** | (0.38) | **0.42** |
| 4.0 | 0.005 | **0.11** | (0.15) | **0.10** | **0.50** | (0.51) | **0.52** |

Here we have used $(d\widehat{\boldsymbol{\ell\mathrm{fdr}}}/\mathbf{dy})\widehat{\boldsymbol{\nu}} = \mathbf{0}$ by homogeneity. As discussed below, formula (5.9) also has some application to the situation where the $z$ values are correlated. □

NOTE. $\mathbf{y}$ is an approximation to the order statistic of the $z$ values, exactly the order statistic if we let bin width $\Delta \to 0$. False discovery rates only depend upon the order statistic, facilitating compact formulas like (5.9).

A formula similar to (5.11) exists for the tail area false discovery rates $\widehat{\ell\mathrm{Fdr}}_k = \log \widehat{\mathrm{Fdr}}(x_k)$,

$$(5.12) \qquad \widehat{\mathrm{cov}}(\widehat{\boldsymbol{\ell\mathrm{Fdr}}}) = B\widehat{G}^{-1}B',$$

$$(5.13) \qquad B = \widehat{S}_0 X_0 \widetilde{G}_0^{-1} \widetilde{X}_0' \widetilde{X} - \widehat{S}X,$$

where, for the case of left-tail $\widehat{\mathrm{Fdr}}$'s, $\widehat{S}_0$ and $\widehat{S}$ are lower triangular matrices,

$$(5.14) \qquad \widehat{S}_{k\ell} = \frac{\widehat{f}_\ell}{\widehat{F}_k} \quad \text{and} \quad \widehat{S}_{0k\ell} = \frac{\widehat{f}_{0\ell}}{\widehat{F}_{0k}} \qquad \text{for } \ell \leq k.$$

Simulation (3.9) for Table 1 was extended to assess the covariance formula (5.9). Table 4 compares the observed standard deviations of $\log \widehat{\mathrm{fdr}}(z)$, now from 250 trials, with the average estimates $\widehat{sd}$ obtained from the square root of the diagonal elements of (5.9). The formula is quite accurate, especially in the empirical null situation; $\widehat{sd}$ was reasonably stable from trial to trial, with coefficient of variation less than 10% for $2.5 \leq z \leq 3.5$.

The "fdr" column in Table 4 is fdr$(z)$, (2.6), based on $f_0(z) \sim N(0,1)$, $f_1(z) \sim N(3,2)$ and $p_0 = 0.9$ as implied by model (3.9). Decisions between null versus nonnull are most difficult in the crucial range $2.5 \leq z \leq 3.5$, where $\mathrm{fdr}(z)$ declines from 0.38 to 0.03. The standard errors for local fdr estimates are about one third again larger than for tail area Fdr, when using the theoretical null. Both give stable estimates in Table 4: a 10% coefficient of variation might mean an estimated $\widehat{\mathrm{fdr}}$ of $0.20 \pm 0.02$, quite tolerable in most large-scale testing situations.

Estimation accuracy is much worse on the empirical null side of the table: a 25% coefficient of variation translates to uncomfortably variable fdr estimates such as $0.20 \pm 0.05$. Now tail area $\widehat{\mathrm{Fdr}}$'s are about one third *more* variable than local $\widehat{\mathrm{fdr}}$'s [and several percent worse still if $F(z)$ in (2.7) is estimated by the usual empirical c.d.f. rather than the parametric estimate corresponding to $\widehat{f}(z)$]. Increasing $N$ by factor $c$ decreases standard errors by roughly $\sqrt{c}$, so taking $N = 6000$ would about halve the boldface values in Table 4. Reducing the degrees of freedom for estimating $f(z)$ from 7 to 5 decreased standard errors by about one third. MLE fitting gave about the same results as central matching here.

"Always use the theoretical null" is not practical advice, even if supplemented by permutation methods. The theoretical or permutation null yields seriously misleading results for the HIV data, as discussed in [10]. Some form of empirical null estimation seems inevitable here, whether using tail area or local false discovery rates (or, for that matter, other simultaneous testing techniques). Of course one should strive for the most efficient possible estimation method, and MLE fitting seems to offer some advantages in this regard.

The equivalent of Lemma 1 when using MLE fitting is derived from the likelihood (4.12). Some definitions in addition to (4.9) are needed:

$$
\begin{aligned}
[a, b] &= \left[\frac{-x_0 - \delta_0}{\sigma_0}, \frac{x_0 - \delta_0}{\sigma_0}\right], \\
H_p(\delta_0, \sigma_0) &= \int_a^b z^p \varphi(z)\,dz, \qquad p = 0, 1, 2, 3, 4, \\
E_p(\delta_0, \sigma_0) &= \frac{\sigma_0^p}{H_0}\bigg[H_p + p\frac{\delta_0}{\sigma_0}H_{p-1} \\
&\qquad + \binom{p}{2}\left(\frac{\delta_0}{\sigma_0}\right)^2 H_{p-2} + \cdots + \left(\frac{\delta_0}{\sigma_0}\right)^p H_0\bigg]
\end{aligned}
$$
(5.15)

[where $H_p = H_p(\delta_0, \sigma_0)$, etc.], and likewise $\widehat{a}, \widehat{b}, \widehat{H}_p, \widehat{E}_p$ for these quantities when $(\delta_0, \sigma_0) = (\widehat{\delta}_0, \widehat{\sigma}_0)$.

Conditional on $N_0$, the number of $z_i$ values observed in $[-x_0, x_0]$, the second factor in (4.12) is a two-parameter exponential family with bivariate

sufficient statistics

$$\begin{pmatrix} Y_1 \\ Y_2 \end{pmatrix} = \frac{1}{N_0} \begin{pmatrix} \sum_{\mathcal{I}_0} z_i \\ \sum_{\mathcal{I}_0} z_i^2 \end{pmatrix}. \tag{5.16}$$

$Y$ has expectation $(E_1(\delta_0, \sigma_0), E_2(\delta_0, \sigma_0))$ and covariance matrix

$$V = \frac{1}{N_0} \begin{pmatrix} E_2 - E_1^2 & E_3 - E_1 E_2 \\ E_3 - E_1 E_2 & E_4 - E_2^2 \end{pmatrix}. \tag{5.17}$$

By definition, an estimate of $f_0^+(z)$ using MLE fitting depends only on the counts "$\mathbf{y}_0$" within $[-x_0, x_0]$, corresponding say to index set $\mathcal{K}_0$, length $K_0$. Let $M_0$ be the $3 \times K_0$ matrix whose $k$th column equals $(1, x_k - Y_1, x_k^2 - Y_2)'$ for $k \in \mathcal{K}_0$. Straightforward but lengthy exponential family calculations produce the influence function of $\widehat{\boldsymbol{\ell}}_0^+ = (\log \widehat{f}_0^+(x_k))$ with respect to $\mathbf{y}_0$, a $K \times K_0$ matrix,

$$\frac{d\widehat{\boldsymbol{\ell}}_0^+}{d\mathbf{y}_0} = \frac{1}{N\widehat{H}_0 \widehat{p}_0} \left[ \mathbf{1}_K \cdot \left(1 - \frac{N_0}{N}\right), \frac{\widehat{U}\widehat{J}\widehat{V}^{-1}}{\widehat{\sigma}_0} \right] M_0, \tag{5.18}$$

where $\mathbf{1}_K$ is a vector of $K$ 1's, $\widehat{V}$ the estimated version of (5.17),

$$\widehat{J} = \widehat{\sigma}_0^2 \begin{pmatrix} 1 & 2\delta_0^2 \\ 0 & \widehat{\sigma} \end{pmatrix}, \tag{5.19}$$

and $U$ the $K \times 2$ matrix with $k$th row

$$u_k = \left( \frac{x_k - \widehat{\delta}_0}{\widehat{\sigma}_0} - \frac{\widehat{H}_1}{\widehat{H}_0}, \frac{(x_k - \widehat{\delta}_0)^2 - \widehat{\sigma}_0^2}{\widehat{\sigma}_0^2} - \frac{\widehat{H}_2 - \widehat{H}_0}{\widehat{H}_0} \right). \tag{5.20}$$

Since $d\widehat{\boldsymbol{\ell}}/d\mathbf{y} = X\widehat{G}^{-1}X'$ as before, (5.9) and (5.18) combine to give $d\widehat{\boldsymbol{\ell}\mathbf{fdr}}/d\mathbf{y}$ for MLE fitting:

LEMMA 2. *The influence function of $\log \widehat{\mathbf{fdr}}$ with respect to $\mathbf{y}$, using MLE fitting, is the $K \times K$ matrix*

$$\frac{d\widehat{\boldsymbol{\ell}\mathbf{fdr}}}{d\mathbf{y}} = \frac{1}{N\widehat{H}_0 \widehat{p}_0} \left[ \mathbf{1} \cdot \left(1 - \frac{N_0}{N}\right), \frac{\widehat{U}\widehat{J}\widehat{V}^{-1}}{\widehat{\sigma}_0} \right] M - X\widehat{G}^{-1}X', \tag{5.21}$$

*where $M$ is the $3 \times K$ matrix with $k$th column $(1, x_k - Y_1, x_k^2 - Y_2)'$ for $k \in \mathcal{K}_0$, and $(0, 0, 0)'$ for $k \notin \mathcal{K}_0$.*

Delta-method estimates of $\mathrm{cov}(\widehat{\ell\mathbf{fdr}})$ for MLE fitting are obtained from Lemma 2 as in (5.11), though the formula does not collapse neatly as in (5.9). We can employ Lemmas 1 and 2 to compare central matching with MLE fitting for the sensitivity of $\widehat{\mathrm{Fdr}}(z)$ to changes in the count vector $\mathbf{y}$. The left panel of Figure 6 compares the two influence functions $d\log \widehat{\mathrm{fdr}}(z = 3.45)/dy_k$, plotted versus bin centerpoint $x_k, k=1,2,\ldots,K$, for the prostate data [$z=3.45$ has $\widehat{\mathrm{fdr}}(z) \doteq 0.2$ for both empirical methods, rather than $\widehat{\mathrm{fdr}}(3.37) = 0.2$ using the theoretical null]. MLE fitting used $x_0 = 2$ in (4.8), accounting for the discontinuities there in its influence curve. The right panel shows the "variance spectrum"

$$(5.22) \qquad S_k(z) = \left(\frac{d\log \widehat{\mathrm{fdr}}(z)}{dy_k}\right)^2 \widehat{\nu}_k, \qquad k=1,2,\ldots,K,$$

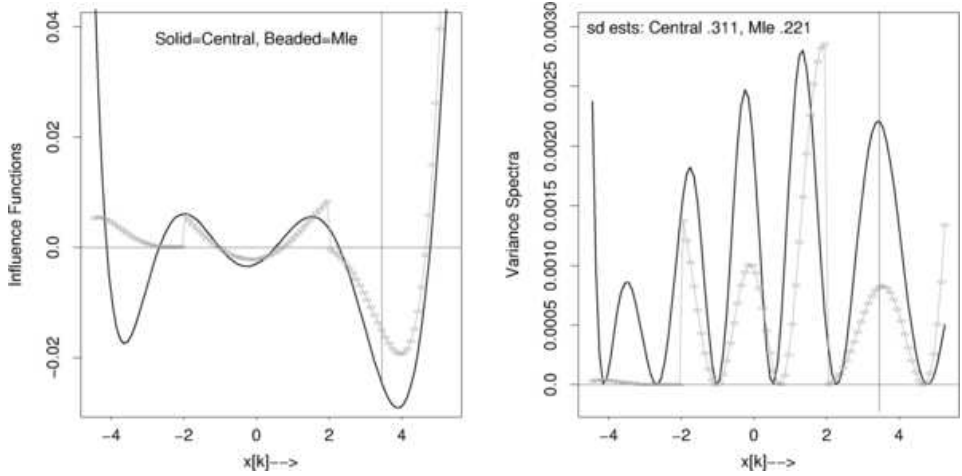

FIG. 6. *Left panel*: *Influence curves* $d\log \widehat{\mathrm{fdr}}(3.45)/dy_k$ *for prostate data, central matching* (solid) *and MLE fitting* (beaded); *plotted versus bin centerpoints* $x_k$. *Right panel*: *variance spectra* (5.22); *MLE has less area, smaller stdev estimate.*

$\widehat{\nu}_k$ the estimated expectation of $y_k$, (5.2). The delta-method estimate of the standard deviation for $\log \widehat{\mathrm{fdr}}(z)$ is

$$(5.23) \qquad \widehat{sd}(z) = \left(\sum_{k=1}^{K} S_k(z)\right)^{1/2}$$

as in (5.11), so variance is proportional to the area under the curve. In this case, MLE fitting has less area and smaller estimated standard deviation. If we were to use an empirical null here, rather than the theoretical null of Figure 3, this would argue for MLE fitting.

Let $\widehat{\xi} = (\widehat{p}_0, \widehat{\delta}_0, \widehat{\sigma}_0)$. Results similar to Lemma 1 and Theorem 1 yield closed-form expressions for the delta-method estimate of $\mathrm{cov}(\widehat{\xi})$. For central matching,

$$(5.24) \qquad \widehat{\mathrm{cov}}(\widehat{\xi}) = D\widehat{G}^{-1}D' - E,$$

$E$ a $3 \times 3$ matrix with $E_{11} = 1/N$ and all other entries 0, and

$$(5.25) \qquad D = \begin{pmatrix} 1 & \widehat{\delta} & \widehat{\sigma}^2 + \widehat{\delta}^2 \\ 0 & \widehat{\sigma}^2 & 2\widehat{\delta}\widehat{\sigma}^2 \\ 0 & 0 & \widehat{\sigma}^3 \end{pmatrix} \widetilde{G}_0^{-1} \widetilde{X}_0' \widetilde{X},$$

$\widehat{G}, \widetilde{G}_0, \widetilde{X}_0$ and $\widetilde{X}$ as in (5.1), (5.3).

The corresponding estimate of $\mathrm{cov}(\widehat{\xi})$ using MLE fitting is

$$(5.26) \qquad \widehat{\mathrm{cov}}(\widehat{\xi}) = aba',$$

with $a$ and $b$ both $3 \times 3$ matrices,

$$(5.27) \qquad a = \begin{pmatrix} 1/\widehat{H}_0, & \mathbf{c}' \\ \mathbf{0} & I_2 \end{pmatrix}, \qquad \mathbf{c}' = -\frac{\widehat{p}_0}{\widehat{\sigma}_0} \left( \frac{\widehat{H}_1}{\widehat{H}_0}, \frac{\widehat{H}_2 - \widehat{H}_0}{\widehat{H}_0} \right),$$

and

$$(5.28) \qquad b = \begin{pmatrix} \widehat{p}_0 \widehat{H}_0 (1 - \widehat{p}_0 \widehat{H}_0 / N), & \mathbf{0}' \\ \mathbf{0}, & \widehat{J}\widehat{V}^{-1}\widehat{J}/N_0 \end{pmatrix},$$

$\widehat{J}$ and $\widehat{V}$ as in (5.17), (5.19). *Locfdr* returns the standard deviation estimates of $\widehat{p}_0, \widehat{\delta}_0$ and $\widehat{\sigma}_0$ based on (5.24) and (5.26).

Model (3.9) presumes that the null genes are *exactly* null. Figure 7 is based on a more relaxed model:

$$(5.29) \quad z_i \overset{\mathrm{ind}}{\sim} N(\mu_i, 1) \quad \text{with } \begin{cases} \mu_i \sim N(0, 0.5^2), & \text{probability } 0.90, \\ \mu_i \sim N(3, 1), & \text{probability } 0.10, \end{cases}$$

$i = 1, 2, \ldots, N = 1500$. In an observational study this might reflect unobserved covariates that jitter even the null cases, as in Section 4 of [8]. In terms of the two-class model (2.2), (5.29) amounts to $p_0 = 0.90$,

$$(5.30) \qquad f_0(z) \sim N(0, 1.12^2) \quad \text{and} \quad f_1(z) \sim N(3, 2).$$

Using the theoretical $N(0,1)$ null in situation (5.29) gives misleading results: $\widehat{\mathrm{fdr}}(z)$ tends to be far too liberal in diagnosing nonnull genes, as shown by the beaded curve in Figure 7. Empirical null estimation gives $\widehat{\mathrm{fdr}}(z)$ estimates much closer to the true curve $\mathrm{fdr}(z) = p_0 f_0(z)/f(z)$ from (5.30). This just says the obvious: empirical methodology correctly estimates $f_0(z)$ in (5.30) [central matching gave $(\widehat{\delta}_0, \widehat{\sigma}_0)$ estimates averaging $(0.02, 1.14)$], which is the whole point of using empirical nulls. Section 4 of [8] discusses what "the correct null" means in this situation, and why it cannot be found by the usual permutation methods.

Our covariance estimates, such as (5.9), were derived assuming independence among the components of $\mathbf{z} = (z, z_2, \ldots, z_N)$ (almost true for the prostate data but not the HIV data). However the influence function formulas have a wider range of applicability. The delta-method estimate of covariance

$$(5.31) \qquad \widehat{\mathrm{cov}}(\widehat{\boldsymbol{\ell}\mathbf{fdr}}) = (d\widehat{\boldsymbol{\ell}\mathbf{fdr}}/d\mathbf{y})\widehat{\mathrm{cov}}(\mathbf{y})(d\widehat{\boldsymbol{\ell}\mathbf{fdr}}/d\mathbf{y})'$$

applies just as well to correlated $z_i$'s. What changes is that (5.10) no longer represents $\mathrm{cov}(\mathbf{y})$.

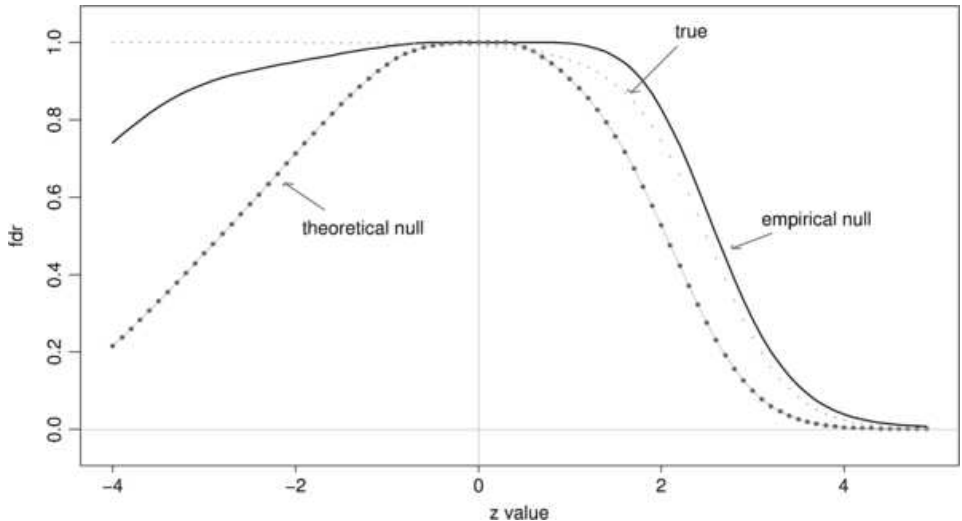

Fig. 7. *Average local false discovery rates* $\widehat{\mathrm{fdr}}(z)$, 250 *replications of* (5.29): *heavy curve using empirical null; beaded curve theoretical null; dots true fdr,* (2.6) *and* (5.30). *Using the theoretical null yields far too many nonnull genes, including some for* $z_i < 0$.

The development in Section 3 of [10] suggests that the estimate

$$(5.32) \qquad \widehat{\mathrm{cov}}(\mathbf{y}) = \mathrm{diag}(\widehat{\boldsymbol{\nu}}),$$

with $\widehat{\nu}_k = N\Delta \widehat{f}(x_k)$ as in (5.2), is still appropriate *in a conditional sense* for the correlated case. Speaking broadly, employing an empirical null amounts to conditioning the estimate $\widehat{\mathrm{fdr}}(z)$ on an approximate ancillary statistic ("A" in [10]), after which (5.31) and (5.32) gives the appropriate conditional covariance. This amounts to using (5.9), or its equivalent for MLE fitting, as stated. More careful estimates of $\widehat{\mathrm{cov}}(\mathbf{y})$ in (5.3) are available in the correlated **z** situation, but the formulas of this section are at least roughly applicable, especially for comparing different estimation techniques.

**6. Summary.** Large-scale simultaneous testing situations, with thousands of hypothesis tests to perform at the same time, are illustrated by the two microarray studies of Figure 1. False discovery rate methods facilitate both size and power calculations, as discussed in Sections 2 and 3, bringing empirical Bayes ideas to bear on simultaneous inference problems. Two types of false discovery rate statistics are analyzed, the more familiar tail area Fdr's introduced by Benjamini and Hochberg [3], and local fdr's, which are better suited for Bayesian interpretation. Power diagnostics may show, as in our examples, that a majority of the nonnull cases cannot be reported as "interesting" to the investigators without including an unacceptably high proportion of null cases.

Fdr methods, either local or tail area, are easy to apply when the appropriate null distribution is known to the statistician from theoretical or permutation considerations. However, it may be clear that the theoretical/permutation null is incorrect, as with the second histogram of Figure 1. Section 4 gives four reasons why this might happen, especially in observational studies. Two methods of estimating an "empirical null" distribution are presented, and formulas for their accuracy derived in Section 5. Using an empirical null decreases the accuracy of false discovery rate methods, both local and tail area, but is unavoidable in situations like the second microarray study. Software in R, *locfdr*, is available through CRAN for carrying out all the fdr size and power calculations.

Department of Statistics
Sequoia Hall
Stanford University
Stanford, California 94305-4065
USA
E-mail: brad@stat.stanford.edu